\documentclass{article}%
\usepackage{graphicx}
\usepackage{amsmath}
\usepackage{amsfonts}
\usepackage{amssymb}%
\newtheorem{theorem}{Theorem}
\newtheorem{acknowledgement}[theorem]{Acknowledgement}

\newtheorem{lemma}[theorem]{Lemma}

\newtheorem{proposition}[theorem]{Proposition}
\newtheorem{remark}[theorem]{Remark}

\newenvironment{proof}[1][Proof]{\textbf{#1.} }{\ \rule{0.5em}{0.5em}}
\begin{document}

\title{The\ time-dependent maximum principle\ for\ systems\ of\ parabolic\ equations
subject to an avoidance set}
\author{Bennett Chow\\University of California, San Diego
\and Peng\ Lu\\McMaster University}
\date{December 22, 2001 (Revised May 17, 2002)}
\maketitle

\section{Introduction}

For \textit{scalar} parabolic equations, maximum principles are well-known
\cite{PW} and have been applied in numerous settings in partial differential
equations and geometric analysis. In the case of \textit{systems }of parabolic
equations, maximum principles are not as well-known and appear to be much less
frequent. Notable exceptions are given by Richard Hamilton \cite{H1},
\cite{H2} and Joel Smoller \cite{S} (see chapter 14)\footnote{We would like to
thank Yung-Sze Choi for informing us of this reference.}. Hamilton's maximum
principle holds for solutions of reaction-diffusion equations (PDE) which are
time-dependent sections of a vector bundle over a Riemannian manifold; in
particular, it holds for the reaction-diffusion equations satisfied by the
curvature operator under Ricci flow. When the convex subsets of the fibers are
independent of time, Hamilton proved such a maximum principle in \cite{H2},
which roughly says that if the convex subsets are preserved by the system of
ordinary differential equations (ODE) associated to the PDE then the convex
subsets are preserved by the PDE. A special case of this result, which applied
to symmetric $2$-tensors, was proved earlier by Hamilton in \cite{H1} and
applied to obtain crucial curvature pinching estimates in his proof that a
compact 3-manifold with positive Ricci curvature converges to a constant
curvature metric under the volume preserving Ricci flow. The general
formulation in \cite{H2} greatly simplified the computations in \cite{H1} and
facilitated the more complicated convex analysis of the ODE associated to the
evolution of the Riemann curvature operator in dimension four in \cite{H5}.

The main purpose of this paper is to prove two extensions of Hamilton's
maximum principle for systems which should be useful for the study of the
Ricci flow and some other geometric evolution equations such as the mean
curvature flow. We present an extension where the convex sets are allowed to
depend on time, we call the extension the time-dependent maximum principle
(see Theorem \ref{Hamilton time-dep weak max thm}). We also present a
souped-up version where both the convex sets are allowed to depend on time and
the convex set may not be preserved by the ODE on a subset of the boundary but
the solution to the PDE avoids that part of the boundary. We call such subsets
of the convex sets, which contain this part of their boundary,
\textit{avoidance sets}, and call this extension\ the time-dependent maximum
principle subject to an avoidance set (see Theorem
\ref{avoidance time-dep weak max thm}).

A special case of this \textit{time-dependent} maximum principle has already
been proved by Hamilton in Theorem 3.3 in section 2.3 of \cite{H5}. In
particular, Hamilton adjoins to the solution $\sigma$ of the PDE the function
$r=\frac{1}{T_{\ast}-t},$ where $T_{\ast}$ is the singularity time, which
trivially satisfies the equation $\frac{\partial r}{\partial t}=\Delta
r+r^{2}$ and applies the time-independent maximum principle to the pair
$\left(  \sigma,r\right)  .$ However, in general this device of adjoining the
function $r=\frac{1}{T-t}$ has the drawback that the sets in space and time to
be preserved may not be convex even though the space slices are. In such a
case Hamilton's proof would not directly apply. In the proof of the somewhat
more general form of the time-dependent maximum principle given in this paper,
we modify Hamilton's original proof of the time-independent maximum principle
in \cite{H2}. The difficulty in this approach is reconciling the
time-dependence of the sets over which one takes the maximum of certain
functions with the framework of Hamilton's \textquotedblleft ODE\ to
PDE\textquotedblright\ formulation. More precisely, when the convex sets
depend on time, the lemma used by Hamilton (see Lemma
\ref{lem taking deriv of sup funct}) on taking the time derivative of the
function $\sup_{s\in\mathcal{S}\left(  t\right)  }g(s,t)$ must have a
correction term, since now the set $\mathcal{S}\left(  t\right)  $ depends on
time, which is difficult to control in the later applications of the lemma.
This is overcome by considering the space-time track of the time-dependent
convex sets and finding suitable splitting of certain quantities which arise
in the study of both the ODE and the PDE (see the proof of Proposition
\ref{prop on ode presev time dep convex set} and the proof of Theorem
\ref{Hamilton time-dep weak max thm}). Formulating the proof this way enables
us to generalize our result to the case where the PDE is subject to an
avoidance set without much difficulty. A special case of this time-dependent
maximum principle has already been applied in \cite{H5} and \cite{H6} to
obtain refined and subtle pointwise curvature estimates. In addition, although
it is not necessary, our time-dependent maximum principle may be used in the
proofs of results in \cite{H3} and \cite{H4}.

Our time-dependent maximum principle subject to avoidance sets is a more
general formulation of a form of the maximum principle implicitly used in the
proof of certain estimates in section 2.3 of \cite{H5} which are used to
detect necks. Some of these estimates, have analogues in dimension three (see
\S 24 of \cite{H4}). There, our souped-up version can also be used to give an
alternate proof of Theorem 24.6 of \cite{H4} (a suggestion of Mao-Pei Tsui).
However, in the study of 4-manifolds with positive isotropic curvature
\cite{H5}, a souped-up version is necessary and is implicitly used by Hamilton.

\begin{acknowledgement}
P. L. would like to thank Pengfei Guan for helpful discussions and Gang Tian
for his constant support and encouragement.
\end{acknowledgement}

\section{Main results}

Let $M^{n}$ be a closed oriented $n$-dimensional manifold with a smooth family
of Riemannian metrics $g(t)$, $t\in\lbrack0,T]$. Let $V\rightarrow M$ be a
real vector bundle with a time-independent bundle metric $h$ and
$\Gamma\left(  V\right)  $ be the vector space of $C^{\infty}$ sections of
$V.$ Let
\[
\nabla(t):\Gamma(V)\rightarrow\Gamma(V\otimes TM^{\ast}),t\in\lbrack0,T]
\]
be a smooth family of time-dependent connections compatible with $h$, that is,%

\[
X[h(\sigma,\tau)]=h(\nabla(t)_{X}\sigma,\tau)+h(\sigma,\nabla(t)_{X}\tau)
\]
for all $X\in TM,\,\sigma,\tau\in\Gamma(V)$ and $t\in\lbrack0,T]$. The
time-dependent Laplacian $\Delta(t)$ acting on a section $\sigma\in\Gamma(V)$
is defined by%

\[
\Delta(t)\sigma=\text{trace}_{g(t)}(\widehat{\nabla}(t)(\nabla(t)\sigma)),
\]
where%
\[
\widehat{\nabla}(t):\Gamma(V\otimes TM^{\ast})\rightarrow\Gamma(V\otimes
TM^{\ast}\otimes TM^{\ast})
\]
is defined using the connection $\nabla(t)$ on $V$ and the Levi-Civita
connection $D(t)$ on $TM^{\ast}$ associated with metric $g(t)$. That is,%

\[
\widehat{\nabla}(t)_{X}(\sigma\otimes\alpha)=(\nabla(t)_{X}\sigma
)\otimes\alpha+\sigma\otimes(D(t)_{X}\alpha)
\]
for all $X\in TM$, $\sigma\in\Gamma(V)$, $\alpha\in\Gamma(TM^{\ast})$.

Let $F:V\times\lbrack0,T]\rightarrow V$ be a fiber preserving map; i.e.,
$F(\sigma,t)$ is a time-dependent vector field defined on the bundle $V$ and
tangent to the fibers. Then we can form a system of reaction-diffusion
equations (PDE)%

\begin{equation}
\frac{\partial}{\partial t}\sigma(x,t)=\Delta(t)\sigma(x,t)+F(\sigma(x,t),t),
\label{simplified heat eq 1.1}%
\end{equation}
where $\sigma(\cdot,t),\,t\in\lbrack0,T]$ are sections of \ $V$. In each fiber
$V_{x}$ the system of ordinary differential equations (ODE) associated to the
PDE (\ref{simplified heat eq 1.1}) obtained by dropping the Laplacian term is%

\begin{equation}
\frac{d}{dt}\sigma_{x}(t)=F(\sigma_{x}(t),t), \label{simplified ode 1.2}%
\end{equation}
where $\sigma_{x}(t)\in V_{x}$.

Let $\mathcal{K}$ be closed subset of $V$. Denote $\mathcal{K}_{x}%
\doteqdot\mathcal{K}\cap V_{x}$. For any initial time $t_{0}\in\lbrack0,T)$ we
say that the solution $\sigma(x,t):t\in\lbrack t_{0},T]$ of the PDE
(\ref{simplified heat eq 1.1}) starts in $\mathcal{K}$ if $\sigma(x,t_{0}%
)\in\mathcal{K}_{x}$ for all $x\in M$. We say that the solution $\sigma(x,t)$
remains in $\mathcal{K}$ for all later times if $\sigma(x,t)\in\mathcal{K}%
_{x}$ for all $x\in M$ and all $t\in(t_{0},T]$. For any $x\in M$ and for any
initial time $t_{0}\in\lbrack0,T)$ we say that the solution $\sigma
_{x}(t):t\in\lbrack t_{0},T]$ of the ODE (\ref{simplified ode 1.2}) starts in
$\mathcal{K}_{x}$ if $\sigma_{x}(t_{0})\in\mathcal{K}_{x}$. We say that the
solution $\sigma_{x}(t)$ remains in $\mathcal{K}_{x}$ for all later times if
$\sigma_{x}(t)\in\mathcal{K}_{x}$ for all $t\in(t_{0},T]$.

One important question is: when will an arbitrary solution of the PDE
(\ref{simplified heat eq 1.1}) which starts in $\mathcal{K}$ at an arbitrary
initial time $t_{0}\in\lbrack0,T)$ remain in $\mathcal{K}$ for all later
times? To answer this question, we need to impose two conditions on
$\mathcal{K}$:\smallskip

\textbf{I}. $\mathcal{K}$ is invariant under parallel translation defined by
the connection $\nabla(t)$ for each $t\in\lbrack0,T]$;

\textbf{II}. In each fiber $V_{x}$ set $\mathcal{K}_{x}$ is closed and
convex.\smallskip

The following theorem is the maximum principle in the time-independent case
(Theorem 4.3 in \cite{H2}).

\begin{theorem}
\label{Hamilton time-indep weak max thm} Let $\mathcal{K}\subset V$ be a
closed subset satisfying conditions \textbf{I} and \textbf{II}. Assume that
$F(\sigma,t)$ is continuous in $t$ and is Lipschitz in $\sigma$. Suppose that
for any $x\in M$ and any initial time $t_{0}\in\lbrack0,T),$\textrm{ }and any
solution $\sigma_{x}(t)$ of the ODE (\ref{simplified ode 1.2}) which starts in
$\mathcal{K}_{x}$ at $t_{0},$ the solution $\sigma_{x}(t)$ will remain in
$\mathcal{K}_{x}$ for all later times. Then for any initial time $t_{0}%
\in\lbrack0,T)$ the solution $\sigma(x,t)$ of the PDE
(\ref{simplified heat eq 1.1}) will remain in $\mathcal{K}$ for all later
times if $\sigma(x,t)$\textrm{ }starts in $\mathcal{K}$ at time $t_{0}$.
\end{theorem}

In applications to the Ricci flow the vector bundle $V$ is a tensor bundle and
the subsets $\mathcal{K}_{x}\subset V_{x},$ which are invariant under the
action of $O\left(  n\right)  ,$ are identified under the isomorphism between
two fibers $V_{x}$ and $V_{y}$ induced by any choice of orthonormal frames in
$TM$ at the two points $x$ and $y.$ The ODEs in $V_{x}$ and $V_{y}$ are also
$O\left(  n\right)  $-invariant and identical under this identification. When
this is the case the requirement on the solutions of the ODE
(\ref{simplified ode 1.2}) in Theorem \ref{Hamilton time-indep weak max thm}
will hold for every fiber if it holds for one fiber.

Next we formulate the maximum principle in the time-dependent case. Let $U$ be
an open subset of $V$ and $\mathcal{K}(t)\subset U$ be a closed subset for
each $t\in\lbrack0,T]$. We impose two conditions on $\mathcal{K}(t)$ for each
$t$:\smallskip

\textbf{III}. $\mathcal{K}(t)$ is invariant under parallel translation defined
by the connection $\nabla(t)$ for each $t\in\lbrack0,T]$;

\textbf{IV}. In each fiber $V_{x}$ set $\mathcal{K}_{x}(t)\doteqdot
\mathcal{K}(t)\cap V_{x}$ is nonempty, closed and convex for each $t\in
\lbrack0,T]$\textrm{.\smallskip}

\noindent

We define the space-time track%

\[
\mathcal{T}\doteqdot\{(v,t)\in V\times\lbrack0,T]:v\in\mathcal{K}%
(t),\text{\ }t\in\lbrack0,T]\},
\]
and define
\[
\mathcal{T}_{x}\doteqdot\mathcal{T}\cap(V_{x}\times\lbrack0,T]).
\]

Let $F:U\times\lbrack0,T]\rightarrow V$ be a fiber preserving map, i.e.,
$F(x,\sigma,t)$ is a time-dependent vector field defined on $U$ and tangent to
the fibers. Let $u(x,t):V_{x}\times T_{x}M^{\ast}\rightarrow V_{x}$ be a
smooth family of bundle maps of diagonal form, i.e.,
\[
u(x,t)(\sigma,dx^{i})=u^{i}(x,t)\cdot\sigma
\]
where $u^{i}(x,t)$ are smooth functions (in applications to the Ricci flow,
$u\equiv0$ ; that is, there is no gradient term.) Then we can form a system of
reaction-diffusion equations (PDE)%

\begin{equation}
\frac{\partial}{\partial t}\sigma(x,t)=\Delta(t)\sigma(x,t)+u(x,t)(\nabla
(t)\sigma(x,t))+F(x,\sigma(x,t),t). \label{general heat eq 1.3}%
\end{equation}
In each fiber $V_{x}$ the associated system of ordinary differential equations
(ODE) is%

\begin{equation}
\frac{d}{dt}\sigma_{x}(t)=F(x,\sigma_{x}(t),t). \label{general ode 1.4}%
\end{equation}

Hamilton's maximum principle is an answer to the following question: For any
$t_{0}\in\lbrack0,T)$ when will the solution $\sigma(x,t),\,t\in\lbrack
t_{0},T]$ of the PDE (\ref{general heat eq 1.3}) which starts in
$\mathcal{K}(t_{0})$, remain in $\mathcal{K}(t)$ for all later times, i.e.,
$\sigma(x,t)\in\mathcal{K}_{x}(t)$ for all $x\in M$ and $\,t\in\lbrack
t_{0},T]$?

In this paper we extend Hamilton's techniques established in \cite{H2} to
the\textit{ time-dependent} case (see Theorem
\ref{Hamilton time-dep weak max thm} below) and use our extension to also
prove a time-dependent maximum principle for the PDE
(\ref{general heat eq 1.3}) subject to an \textit{avoidance set} (see Theorem
\ref{avoidance time-dep weak max thm} below).

\begin{theorem}
\label{Hamilton time-dep weak max thm} Let $\mathcal{K}(t)\subset
V,t\in\lbrack0,T]$ be closed subsets which satisfy conditions\textrm{
}\textbf{III}\textrm{ }and\textrm{ }\textbf{IV}\textrm{ }above, and such that
the space-time track $\mathcal{T}$\textrm{ }is\textrm{ }closed. Assume that
$u(x,t):V_{x}\times T_{x}M^{\ast}\rightarrow V_{x}$ is a smooth family of
bundle maps of diagonal form and assume that $F(x,\sigma,t)$ is continuous in
$x,\,t$ and is Lipschitz in $\sigma$. Suppose that, for any $x\in M$ and any
initial time $t_{0}\in\lbrack0,T),$ any solution $\sigma_{x}(t)$ of the ODE
(\ref{general ode 1.4}) which starts in $\mathcal{K}_{x}(t_{0})$ will remain
in $\mathcal{K}_{x}(t)$ for all later times, i.e., $\sigma_{x}(t)\in
\mathcal{K}_{x}(t)$ for all $t\in\lbrack t_{0},T]$. Then for any initial time
$t_{0}\in\lbrack0,T)$ the solution $\sigma(x,t):t\in\lbrack t_{0},T]$ of the
PDE (\ref{general heat eq 1.3}) will remain in $\mathcal{K}(t)$ for all later
times if $\sigma(x,t)$ starts in $\mathcal{K}(t_{0})$ at time $t_{0}.$
\end{theorem}

Theorem \ref{Hamilton time-dep weak max thm} is the maximum principle in
time-dependent case (the convex set $\mathcal{K}(t)$ depends on time $t$).
Special cases of this result have been proved by Hamilton and applied to the
study of the Ricci flow (see for example section 2.2 and section 2.3 in
\cite{H5}).

There is also a souped-up version of the maximum principle for systems of
reaction-diffusion equations in time-dependent case. The idea is that for
applications sometimes we are in the situation where the reason that a convex
set $\mathcal{K}(t)$ is not preserved by the ODE (\ref{general ode 1.4}) is
that the solution wants to escape from a certain part of the convex set (which
we will call the avoidance set $\mathcal{A}(t)\subset\mathcal{K}(t)$). In this
case, if we assume that any solution $\sigma(x,t):t\in\lbrack t_{0},T]$ to the
PDE (\ref{general heat eq 1.3}) starts in $\mathcal{K}(t_{0})\backslash
\mathcal{A}(t_{0})$ and assume that the solution $\sigma(x,t)$ does not enter
subsets $\mathcal{A}(t)$ for all $t\geq t_{0}$ (i.e., $\sigma(x,t)\notin
\mathcal{A}_{x}(t)\doteqdot\mathcal{A}(t)\cap V_{x}$ for all $x\in M$ and all
$t\geq t_{0}$\ ), then\textrm{ }$\sigma(x,t)$ remains in $\mathcal{K}(t)$ for
$t\geq t_{0}$. A typical example where this happens is when the solution to
the Ricci flow is assumed not to have any necklike points (see Theorem 3.3 and
3.4 in section 2.3 of \cite{H5}).

We define the avoidance space-time track
\[
\mathcal{AT}\ \doteqdot\left\{  (v,t)\in V\times\lbrack0,T]:v\in
\mathcal{A}(t),\;t\in\lbrack0,T]\right\}  ,
\]
and define%
\[
\mathcal{AT}_{x}\doteqdot(\mathcal{AT)}\cap(V_{x}\times\lbrack0,T]).
\]

\begin{theorem}
\label{avoidance time-dep weak max thm} Let $\mathcal{K}(t)\subset
V,t\in\lbrack0,T]$ be closed subsets which satisfy conditions\textrm{
}\textbf{III}\textrm{ }and\textrm{ }\textbf{IV}\textrm{ }above, and such that
the space-time and the avoidance space-time tracks $\mathcal{T}$\textrm{ }and
$\mathcal{AT}$\ are\textrm{ }closed. Assume that $u(x,t):V_{x}\times
T_{x}M^{\ast}\rightarrow V_{x}$ is a smooth family of bundle maps of the
diagonal form and that $F(x,\sigma,t)$ is continuous in $x,t$ and is Lipschitz
in $\sigma.$\textrm{ }Suppose that for any $x\in M,\,t_{0}\in\lbrack0,T)$ and
any solution $\sigma_{x}(t)$ of the ODE (\ref{general ode 1.4}) with initial
condition $\sigma_{x}(t_{0})\in\mathcal{K}_{x}(t_{0})\backslash\mathcal{A}%
_{x}(t_{0})$, either $\sigma_{x}(t)\in\mathcal{K}_{x}(t)$ for all $t\geq
t_{0}$, or there is time $t_{1}$ such that $\sigma_{x}(t)\in\mathcal{K}%
_{x}(t)\backslash\mathcal{A}_{x}(t)$ for all $t_{0}\leq t\leq t_{1}$ and
$\sigma_{x}(t_{1})\in\mathcal{A}_{x}(t_{1})$. Then for any $t_{0}\in
\lbrack0,T)$ and any solution $\sigma(x,t):t\in\lbrack t_{0},T]$ of the PDE
(\ref{general heat eq 1.3}) satisfying initial condition $\sigma(x,t_{0}%
)\in\mathcal{K}_{x}(t_{0})\backslash\mathcal{A}_{x}(t_{0})$ for all $x\in
M$\ and satisfying $\sigma(x,t)\notin\mathcal{A}_{x}(t)$ for all $x\in M$ and
all $t\geq t_{0},$ the solution\ $\sigma(x,t)$\ will remain in $\mathcal{K}%
(t)$ for later times.
\end{theorem}

\section{Hamilton's proof of Theorem \ref{Hamilton time-indep weak max thm}%
\label{thm1,1 section}}

In order to present the proof of our main result, Theorem
\ref{Hamilton time-dep weak max thm}, more clearly and to exhibit the
differences between the proofs of Theorem \ref{Hamilton time-dep weak max thm}
and Hamilton's Theorem \ref{Hamilton time-indep weak max thm}, we will first
review his proof of Theorem \ref{Hamilton time-indep weak max thm} in this
section. This will enable us to omit the common parts of the two proofs when
we prove Theorem \ref{Hamilton time-dep weak max thm}.

\begin{theorem}
\label{modified Hamilton time-indep weak max thm} Let $\mathcal{K}\subset V$
be a closed subset\textrm{ }satisfying conditions \textbf{I} and \textbf{II}.
Assume that $F(x,\sigma,t)$ is continuous in $x,t$ and is Lipschitz in
$\sigma$. Suppose that for any $x\in M$ and any initial time $t_{0}\in
\lbrack0,T),$ and any solution $\sigma_{x}(t)$ of the ODE
(\ref{general ode 1.4}) which starts in $\mathcal{K}_{x}$ at $t_{0},$ the
solution $\sigma_{x}(t)$ will remain in $\mathcal{K}_{x}$ for all later times.
Then for any initial time $t_{0}\in\lbrack0,T)$ the solution $\sigma(x,t)$ of
the PDE (\ref{general heat eq 1.3}) will remain in $\mathcal{K}$ for all later
times if $\sigma(x,t)$ starts in $\mathcal{K}$ at time $t_{0}$.
\end{theorem}

\begin{remark}
The above result is slightly more general than Theorem
\ref{Hamilton time-indep weak max thm} in that it allows for a gradient term
in the equation (along the lines of the maximum principle for symmetric
$2$-tensors in \cite{H1}). This does not affect Hamilton's proof in \cite{H2}.
\end{remark}

Before proving Theorem \ref{modified Hamilton time-indep weak max thm}, we
need to recall three lemmas, which are essentially in \cite{H2}.

Let $f:[a,b]\rightarrow\mathbb{R}$ be a function. Then we define $\frac
{d^{+}f(t)}{dt}$ at $t\in\lbrack a,b)$ to be the lim sup of forward difference quotients:%

\begin{equation}
\frac{d^{+}f(t)}{dt}=\limsup_{s\rightarrow0^{+}}\,\frac{f(t+s)-f(t)}{s}.
\label{eq 2.1 def of forward derivative}%
\end{equation}

\begin{lemma}
\label{lem on forward deriv} Suppose function $f:[a,b]\rightarrow\mathbb{R}$
is left lower semi-continuous and right-continuous with $f(a)\leq0$. Assume either

\begin{enumerate}
\item[(i)] $\frac{d^{+}f(t)}{dt}\leq0$ when $f(t)\geq0$ on $(a,b)$, or

\item[(ii)] for some constant $C<+\infty,$ $\frac{d^{+}f(t)}{dt}\leq C\cdot
f(t)$ when $f(t)\geq0$ on $(a,b)$.\newline\noindent Then $f(t)\leq0$ on
$[a,b]$.
\end{enumerate}
\end{lemma}

\begin{proof}
By checking the proof of Lemma 3.1 in \cite{H2}, one can prove:

\begin{description}
\item[ Sublemma.] Suppose function $f(t)\mathrm{:}[a,b]\rightarrow\mathbb{R}$
is left lower semi-continuous and right-continuous with $f(a)\leq0$. Assume
$\frac{d^{+}f(t)}{dt}\leq0$ when $f(t)\geq0$ on $[a,b)$. Then $f(t)\leq0$ on
$[a,b]$.
\end{description}

\noindent Hypothesis (i) in Lemma \ref{lem on forward deriv} is a little
weaker than the hypothesis in the sublemma since we do not require the
inequality to hold at the left endpoint $a$. By the right continuity at $t=a,$
given any $\varepsilon>0,$ there exists $\delta>0$ such that $f(t)\leq
\varepsilon$ on $[a,a+\delta].$ Now $f\left(  a+\delta\right)  -\varepsilon
\leq0$ and $\frac{d^{+}\left[  f(t)-\varepsilon\right]  }{dt}\leq0$ when
$f(t)-\varepsilon\geq0$ for $t\in\lbrack a+\delta,b).$ Hence we may apply the
sublemma to $f(t)-\varepsilon$ to obtain $f(t)-\varepsilon\leq0$ on
$[a+\delta,b]$. We conclude that for any $\varepsilon>0,$ we have $f\left(
t\right)  \leq\varepsilon$ for all $t\in\left[  a,b\right]  .$ This proves the
lemma under hypothesis (i).

To prove the lemma under hypothesis (ii), we set $g(t)=e^{-C\cdot t}\cdot
f(t)$. Then $\frac{d^{+}g(t)}{dt}\leq0$ when $g(t)\geq0$ on $(a,b)$. Applying
the lemma under hypothesis (i) to $g(t)$, we get $g\left(  t\right)  \leq0$ on
$\left[  a,b\right]  .$ This implies $f\left(  t\right)  \leq0$ on $\left[
a,b\right]  .\medskip$
\end{proof}

The second lemma below gives a useful characterization of when systems of
ordinary differential equations preserve closed convex sets in Euclidean
space. Let $\mathcal{J}\subset\mathbb{R}^{n}$ be a closed convex subset and
$\partial\mathcal{J}$ be the boundary of $\mathcal{J}$ in $\mathbb{R}^{n}$.
For any $v\in\partial\mathcal{J}$ we define the tangent cone $C_{v}%
\mathcal{J}$ of $\mathcal{J}$ at $v$ to be the smallest convex cone in
$\mathbb{R}^{n}$ with vertex at $v$ which contains $\mathcal{J}$.

\begin{lemma}
\label{lem on ode presev convex set} Let $U\subset\mathbb{R}^{n}$ be an open
subset and $\mathcal{J}\subset U$ be a closed convex subset. Consider the ODE%
\begin{equation}
\frac{d\tau}{dt}=F(\tau,t), \label{general ode 1.5 in Euclidean}%
\end{equation}
where $F:U\times\lbrack0,T]\rightarrow\mathbb{R}^{n}$ is continuous\ in $t$
and is Lipschitz in $\tau$. Then the following two statements are equivalent.

(i) For any initial time $t_{0}\in\lbrack0,T)$, any solution of the ODE
(\ref{general ode 1.5 in Euclidean}) which starts in $\mathcal{J}$ at $t_{0}$
will remain in $\mathcal{J}$ for all later times;

(ii) $v+F(v,t)\in C_{v}\mathcal{J}$ for all $v\in\partial\mathcal{J}$ and
$t\in\lbrack0,T)$.
\end{lemma}

\begin{proof}
This is Lemma 4.1 in \cite{H2}. The fact that $F(\tau,t)$ depends on time $t$
does not pose any difficulties for the original proof.\medskip
\end{proof}

The third lemma gives a general principle on how to take the derivative of a
$sup$-function which plays an important role in proving Theorem
\ref{modified Hamilton time-indep weak max thm}. Note that $\mathcal{S}$ in
the lemma below is independent of time $t$.

\begin{lemma}
\label{lem taking deriv of sup funct} Let $\mathcal{S}$ be a sequentially
compact topological space and let $g:\mathcal{S}\times\lbrack a,b]\rightarrow
\mathbb{R}$ be a function. If $g$ is continuous in $s$ and $t$ and
$\frac{\partial g}{\partial t}$ is continuous in $s$ and $t$, then the
function $f:[a,b]\rightarrow\mathbb{R}$ defined by%
\[
f(t)=\sup_{s\in\mathcal{S}}g(s,t),
\]
is Lipschitz and%
\[
\frac{d^{+}f(t)}{dt}\leq\sup\{\frac{\partial g}{\partial t}(s,t):s\in
\mathcal{S}\text{ satisfies }g(s,t)=f(t)\}.
\]

\end{lemma}

\begin{proof}
See Lemma 3.5 in \cite{H2}. \medskip
\end{proof}

The rest of this section will be devoted to proving Theorem
\ref{modified Hamilton time-indep weak max thm}. As remarked in the proof of
Lemma 4.1 on p. 160 of \cite{H2} we may assume that $\mathcal{K}$ is compact.
For if there were a counterexample $\sigma_{0}(x,t)$ for $t\in\lbrack
t_{0},T]$, then $\sigma_{0}(x,t)$ will be contained in $V(r)$ for some $r$
large enough, where $V(r)$ is the tubular neighborhood of the zero section in
$V$ whose intersection with each fiber $V_{x}$ is a ball of radius $r$ around
origin measured by metric $h.$ Let $\eta$ be a cut-off function on $V$ which
equals 1 on $V(r)$ and equals to zero on $V\backslash V(2r)$. Then we can
modify the PDE (\ref{general heat eq 1.3}) as%

\begin{equation}
\frac{\partial}{\partial t}\sigma(x,t)=\Delta(t)\sigma(x,t)+u(x,t)(\nabla
(t)\sigma(x,t))+\eta(\sigma(x,t))\cdot F(x,\sigma(x,t),t).
\label{cutoff of general heat eq 1.3}%
\end{equation}
Note that the paths of the counterexample solution $\sigma_{0}(x,t)$ do not
change inside $V(r)$, hence $\sigma_{0}(x,t)$ still is a solution of
(\ref{cutoff of general heat eq 1.3}). If we intersect $\mathcal{K}$ with
$V(2r)$, we get a counterexample of Theorem
\ref{modified Hamilton time-indep weak max thm} for
(\ref{cutoff of general heat eq 1.3}) with the closed compact convex set
$V(2r)\cap\mathcal{K}\neq\emptyset$ replacing $\mathcal{K}$, since using Lemma
\ref{lem on ode presev convex set} it is easy to check that the ODE%
\[
\frac{d}{dt}\sigma_{x}(t)=\eta(\sigma_{x}(t))\cdot F(x,\sigma_{x}(t),t)
\]
and $V(2r)\cap\mathcal{K}$ satisfy the assumption of Theorem
\ref{modified Hamilton time-indep weak max thm}.

Now we assume that $\mathcal{K}$ is compact. We define the distance between
$\sigma\in V_{x}$ and $v\in V_{x}$ using the metric $h$ and denote it by
$|\sigma-v|$. We will prove Theorem
\ref{modified Hamilton time-indep weak max thm} by contradiction. Suppose we
have a solution $\sigma(x,t)$ of the PDE (\ref{general heat eq 1.3}) which
starts with $\sigma(x,t_{0})\in\mathcal{K}_{x}$ for all $x\in M$ and which
goes out of $\mathcal{K}$ at some time $t_{2}$. Since $\mathcal{K}$ is closed,
we can find a time $t_{1}\geq t_{0}$ such that $\sigma(x,t_{1})\in
\mathcal{K}_{x}$ for all $x\in M,$ and for any $t\in(t_{1},t_{2})$ there is
$x$ such that $\sigma(x,t)\notin\mathcal{K}_{x}$. Below we will focus on the
time interval $[t_{1},t_{2}]$.

Define the function%

\[
f(t)=\sup_{x\in M}d(\sigma(x,t),\mathcal{K}_{x})=\sup_{x\in M}\inf
_{v\in\mathcal{K}_{x}}|\sigma(x,t)-v|\text{\textrm{ }for }t\in\lbrack
t_{1},t_{2}].
\]
We have $f(t_{1})=0$ and $f(t)>0$ for $t\in(t_{1},t_{2}]$ by assumption. It is
easy to check using condition \textbf{I} that $f(t)$ is a continuous function
of $t$. Below we will prove that there is a constant $C<\infty$ such that
$\frac{d^{+}f(t)}{dt}\leq C\cdot f(t)$ for $t\in(t_{1},t_{2})$. Once this is
proved, then $f(t)\leq0$ for $t\in\lbrack t_{1},t_{2}]$ by Lemma
\ref{lem on forward deriv}(ii). Hence $\sigma(x,t)\in\mathcal{K}_{x}$ for all
$x\in M$ and all $t\in(t_{1},t_{2}]$, we get the required contradiction.

For\ any $v\in\partial\mathcal{K}_{x}$, let $S_{v}\subset V_{x}$ be the set of
outward normal directions $n$ of the supporting hyperplanes of $\mathcal{K}%
_{x}$ at $v$; we require that $n$ be unit with respect to the metric $h$.
Then, since $\mathcal{K}$ is nonempty and for each $t\in(t_{1},t_{2}),$
$\sigma(x,t)$ is not in $\mathcal{K}_{x}$ for some $x\in M$, it is well-known that%

\begin{equation}
f(t)=\sup_{x\in M}\sup_{v\in\partial\mathcal{K}_{x}}\sup_{n\in S_{v}}%
n\cdot(\sigma(x,t)-v), \label{eq 2.3 for sup funct from dist}%
\end{equation}

\noindent where $\cdot$ is the inner product in $V_{x}$ defined by the metric
$h$.\textrm{ }Define the set%

\[
\mathcal{S}=\{(x,v,n):x\in M,v\in\partial\mathcal{K}_{x},n\in S_{v}\}
\]
and the function%
\[
g((x,v,n),t)=n\cdot(\sigma(x,t)-v),
\]
then
\[
f(t)=\sup_{(x,v,n)\in\mathcal{S}}g((x,v,n),t).
\]

Note that $S$ is a compact subset of $V\otimes V$ independent of time $t$, we
can apply Lemma \ref{lem taking deriv of sup funct} and get for any
$t\in(t_{1},t_{2})$%

\[
\frac{d^{+}f(t)}{dt}\leq\sup\frac{\partial}{\partial t}[n\cdot(\sigma
(x,t)-v)],
\]
where the sup is over all $(x,v,n)\in\mathcal{S}$ such that $n\cdot
(\sigma(x,t)-v)=f(t)$; in particular we have $|\sigma(x,t)-v|=f(t)$ for these
$(x,v,n)$. We compute at these $(x,v,n)$
\begin{align*}
&  \frac{\partial}{\partial t}[n\cdot(\sigma(x,t)-v)]\newline=n\cdot
(\frac{\partial}{\partial t}\sigma(x,t))\\
&  =n\cdot\lbrack\Delta(t)\sigma(x,t)]+n\cdot\lbrack u(x,t)(\nabla
(t)\sigma(x,t))]+n\cdot F(x,\sigma(x,t),t).
\end{align*}
By the assumption of Theorem \ref{modified Hamilton time-indep weak max thm}
and Lemma \ref{lem on ode presev convex set} we have $v+F(x,v,t)\in
C_{v}\mathcal{K}_{x}$. Hence $n\cdot F(x,v,t)\leq0$ for any $n\in S_{v}$ and
any $t\in(t_{1},t_{2})$. We have%
\begin{align*}
&  n\cdot F(x,\sigma(x,t),t)\\
&  \leq n\cdot F(x,\sigma(x,t),t)-n\cdot F(x,v,t)\newline\leq n\cdot\lbrack
F(x,\sigma(x,t),t)-F(x,v,t)]\newline\\
&  \leq|F(x,\sigma(x,t),t)-F(x,v,t)|\newline\leq C\cdot|\sigma(x,t)-v|\newline%
=C\cdot f(t),
\end{align*}
where $C$ is some constant from the assumption that $F(x,\sigma,t)$ is
Lipschitz in $\sigma$.

We claim that
\begin{align*}
n\cdot\lbrack u(x,t)(\nabla(t)\sigma(x,t))]  &  =0,\\
n\cdot\lbrack\Delta(t)\sigma(x,t)]  &  \leq0,
\end{align*}
which will be proved in a moment. This shows%

\[
\frac{d^{+} f(t)}{dt} \leq C \cdot f(t) \text{ on } (t_{1},t_{2}).
\]

We are left to prove the claim. We will prove $n\cdot\lbrack u(x,t)(\nabla
(t)\sigma(x,t))]=0$ and $n\cdot\lbrack\Delta(t)\sigma(x,t)]\leq0$ together.
Recall that $(x,v,n)$ satisfies $n\cdot(\sigma(x,t)-v)=f(t)$. If we extend a
vector in the bundle $V$ from a point $x$ by parallel translation along
geodesics emanating radially out of $x$, we get a smooth section of the bundle
on some small neighborhood of $x$ such that all the symmetrized covariant
derivatives of the section at $x$ are zero. Let $y$ be an arbitrary point in
some small neighborhood $U_{x}$ of $x$. We extend $v\in\partial\mathcal{K}%
_{x}$ and $n\in V_{x}$ in this manner using the connection $\nabla(t)$ to get
$v_{y}$ and $n_{y}$. Since the connection $\nabla(t)$ is compatible with the
metric $h$ we continue to have $|n_{y}|=1$, and since $\mathcal{K}$ is
invariant under parallel translation we have $v_{y}\in\partial\mathcal{K}_{y}$
and $n_{y}\in S_{v_{y}}$ for $\mathcal{K}_{y}$ at $v_{y}$. Therefore%

\[
n_{y}\cdot(\sigma(y,t)-v_{y})\leq f(t),
\]
for all $y\in U_{x}$. It follows that function $n_{y}\cdot(\sigma(y,t)-v_{y})$
of $y\in U_{x}$ has a local maximum at $y=x$. So%
\begin{align*}
\frac{\partial}{\partial y^{i}}[n_{y}\cdot(\sigma(y,t)-v_{y})]  &  =0\text{ at
}y=x,\\
\Delta(t)[n_{y}\cdot(\sigma(y,t)-v_{y})]  &  \leq0\text{ at }y=x.\newline%
\end{align*}

Let $\nabla_{t,i}$ be the covariant derivative in direction $\frac{\partial
}{\partial y^{i}}$ defined by the connection $\nabla(t)$. Since $v_{y}$ and
$n_{y}$ have their symmetrized covariant derivatives equal to zero at $y=x$,
so $\nabla_{t,i}n_{y}=\nabla_{t,i}v_{y}=0$ and $\Delta(t)n_{y}=\Delta
(t)v_{y}=0$ at $y=x$. Hence%
\[
n\cdot\lbrack\nabla_{t,i}\sigma(x,t)]=0,\hspace{0.5in}n\cdot\lbrack
\Delta(t)\sigma(x,t)]\leq0.
\]
Then%
\begin{align*}
n\cdot\lbrack u(x,t)(\nabla(t)\sigma(x,t))]  &  =n\cdot\lbrack\sum_{i}%
u^{i}(x,t)\cdot\nabla_{t,i}\sigma(x,t)]\newline\\
&  =\sum_{i}u^{i}(x,t)\cdot(n\cdot\lbrack\nabla_{t,i}\sigma(x,t)])=0.
\end{align*}
The claim is proved and so is Theorem
\ref{modified Hamilton time-indep weak max thm}.

\section{Proof of Theorem \ref{Hamilton time-dep weak max thm}%
\label{sect time-dependent thm}}

Throughout this section we will use the same index notation $i$ to denote a
sequence or its subsequence or the subsequence of its subsequence. Our
arguments below will involve taking subsequences from time to time. The
convention will simplify our notations. Before proving Theorem
\ref{Hamilton time-dep weak max thm}, we first formulate a useful
characterization of when systems of ordinary differential equations preserve
time-dependent closed convex sets in Euclidean space, i.e., a time-dependent
version of Lemma \ref{lem on ode presev convex set}.

Let $\mathcal{J}(t)\subset\mathbb{R}^{n},0\leq t\leq T$ be a family of
non-empty closed convex subsets. Define the space-time track%

\[
\mathcal{L}=\{(v,t)\in\mathbb{R}^{n}\times\mathbb{R}:v\in\mathcal{J}%
(t),\,0\leq t\leq T\}.
\]
For each $(v,t)\in\mathcal{L}$ we define a time-like tangent cone in the
forward direction of $\mathcal{L}$ at $(v,t)$ and denote it by\textrm{
}$C_{(v,t)}\mathcal{L}.$ $C_{(v,t)}\mathcal{L}$ consists of all $(W,1)\in
\mathbb{R}^{n}\times\mathbb{R}$ satisfying the following condition: For any
sequence $s_{i}\rightarrow0^{+}$ (i.e., $s_{i}$ approaches to zero from
positive side), there is a subsequence of $s_{i}$ and vectors $W_{i}%
\rightarrow W$ such that points $(v+s_{i}W_{i})\in\mathcal{J}(t+s_{i})$. Note
that the definition is stronger than the conventional definition where one
sequence of $s_{i}$ is enough. When $\mathcal{J}(t)=\mathcal{J}$ is
independent of time $t,$ then $C_{(v,t)}\mathcal{L}=\{C_{v}\mathcal{J}%
-v\}\times\{1\}.$

\begin{proposition}
\label{prop on ode presev time dep convex set} Let $U\subset R^{n}$\textrm{
}be an open subset and $\mathcal{J}(t)\subset U,0\leq t\leq T$\textrm{ }be a
family of non-empty closed convex subsets such that the space-time track
$\mathcal{L}$ is closed. Consider the ODE%
\begin{equation}
\frac{d\tau}{dt}=F(\tau,t), \label{eq 2.5 for ode in Euclidean space}%
\end{equation}
where $F:U\times\lbrack0,T]\rightarrow\mathbb{R}^{n}$ is continuous in\textrm{
}$t$ and is Lipschitz in $\tau$. Then the following two statements are equivalent.

(i) For any initial time $t_{0}\in\lbrack0,T)$, any solution of the
ODE\textrm{ }(\ref{eq 2.5 for ode in Euclidean space}) which starts in
$\mathcal{J}(t_{0})$ at time $t_{0}$ will remain in $\mathcal{J}(t)$ for all
later times;

(ii) $(F(v,t),1)\in C_{(v,t)}\mathcal{L}$ for all $(v,t)\in\partial
\mathcal{L}$, where $\partial\mathcal{L}$ is the boundary of $\mathcal{L}%
\subset\mathbb{R}^{n+1}$.
\end{proposition}

\textbf{Proof.} $(i)\Rightarrow(ii)$. For any $(v_{0},t_{0})\in\partial
\mathcal{L}$, we consider the solution of
(\ref{eq 2.5 for ode in Euclidean space}) with initial condition $\tau
(t_{0})=v_{0}$. (i) implies that $\tau(t_{0}+s)\in\mathcal{J}(t_{0}+s)$ for
any $s\in\lbrack0,T-t_{0}]$. Hence%
\[
\lim_{s\rightarrow0^{+}}\frac{(\tau(t_{0}+s),t_{0}+s)-(\tau(t_{0}),t_{0})}%
{s}=(F(v_{0},t_{0}),1)\in C_{(v_{0},t_{0})}\mathcal{L}.
\]

$(ii)\Rightarrow(i)$. We prove it by contradiction. We will not assume
$\mathcal{L}$ to be compact. Suppose we have a solution $\tau(t)$ starting
with $\tau(t_{0})\in\mathcal{J}(t_{0})$ and going out of $\mathcal{L}$ at some
time $t_{2}$, i.e., $\tau(t_{2})\notin\mathcal{J}(t_{2})$. Since $\mathcal{L}$
is closed, we can find a time $t_{1}$ such that $\tau(t_{1})\in\mathcal{J}%
(t_{1})$ and $\tau(t)\notin\mathcal{J}(t)$ for all $t\in(t_{1},t_{2}%
)$.\textrm{ }Below we will focus on the time interval $[t_{1},t_{2}]$.

Let $\partial\mathcal{J}(t)$ be the boundary of $\mathcal{J}(t)\subset
\mathbb{R}^{n}$. Define the function%
\[
l(t)=d(\tau(t),\mathcal{J}(t))\text{\textrm{ }for\textrm{ }}t\in\lbrack
t_{1},t_{2}]
\]
where $d$ is the Euclidean distance on $\mathbb{R}^{n}$. It is clear that
$l(t_{1})=0$ and $l(t)>0$ for $t\in(t_{1},t_{2}]$. Because $\mathcal{L}$ is
not assumed to be a domain with smooth boundary, the function $l(t)$ is not
necessarily continuous.

\begin{lemma}
\label{lem on dist is lower cont in Eucl} Let $\mathcal{J}(t)\subset U,0\leq
t\leq T$\textrm{ }be a family of non-empty closed convex subsets. If the
space-time track $\mathcal{L}$\textrm{ }is closed and satisfies (ii) in
Proposition \ref{prop on ode presev time dep convex set}, then $l(t)$ is left
lower semi-continuous and is right continuous on $[t_{1},t_{2}]$.
\end{lemma}

\textbf{Proof} of the lemma. To see that $l(t)$ is lower semi-continuous, for
any $t\in\lbrack t_{1},t_{2}]$ and any $s_{i}\rightarrow0$ with $t+s_{i}%
\in(t_{1},t_{2}]$, we choose $v_{i}\in\partial\mathcal{J}(t+s_{i})$ such
that\textrm{ }
\[
l(t+s_{i})=d(\tau(t+s_{i}),v_{i}).
\]
Then either a subsequence $v_{i}$ will converge to some $v_{\infty}%
\in\mathcal{J}(t)$ since $\mathcal{L}$ is closed, or $v_{i}$ will diverge to
$\infty$. In the case of convergence, we have $l(t+s_{i})\rightarrow
d(\tau(t),v_{\infty})\geq l(t)$. The lower semi-continuity is true. In the
case of divergence, then $l(t+s_{i})\rightarrow+\infty$. Since ${\partial
}\mathcal{J}(t)$ is nonempty, $l(t)$ is finite\textrm{. }Hence the lower
semi-continuity of $l(t)$ is also true.

To prove the right-continuity of $l(t)$, it suffices to prove the upper
right-continuity. We will use (ii) in Proposition
\ref{prop on ode presev time dep convex set} which actually puts some
restriction on the space-time track $\mathcal{L}$. It follows from
$(\tau(t),t)\notin\mathcal{L}$ for $t>t_{1}$ that $(\tau(t_{1}),t_{1}%
)\in\partial\mathcal{L}$. We denote $\tau(t_{1})$ by $v_{t_{1}}$. For any
$t\in(t_{1},t_{2})$ it follows from\textrm{ }$\tau(t)\notin\mathcal{J}(t)$
that there is $v_{t}\in\partial\mathcal{J}(t)$ such that $l(t)=d(\tau
(t),v_{t}).$\textrm{ }Hence for any $t\in\lbrack t_{1},t_{2})$ we can find
$v_{t}\in\mathcal{J}(t)$ such that $l(t)=d(\tau(t),v_{t})$ and $(v_{t}%
,t)\in\partial\mathcal{L}$. By (ii) $(F(v_{t},t),1)\in C_{(v_{t}%
,t)}\mathcal{L}$. If we fix a $t\in\lbrack t_{1},t_{2}),$\textrm{ }then for
any sequence $s_{i}\rightarrow0^{+}$ we can find a subsequence $s_{i}$ such
that $(v_{t}+s_{i}W_{i})\in\mathcal{J}(t+s_{i})$ and $W_{i}\rightarrow
F(v_{t},t)$. So
\[
l(t+s_{i})\leq d(\tau(t+s_{i}),v_{t}+s_{i}W_{i}).
\]
Letting $i\rightarrow\infty$, we get $\lim\sup_{i\rightarrow+\infty}%
l(t+s_{i})\leq d(\tau(t),v_{t})=l(t)$. Hence $\lim\sup_{i\rightarrow+\infty
}l(t+s_{i})=l(t)$ by the lower semi-continuity of $l(\cdot)$. The lemma is proved.

Now we go back to the proof of $(ii)\Rightarrow(i)$ in Proposition
\ref{prop on ode presev time dep convex set}.\ Below we will prove that there
is some constant $C<\infty$\textrm{ }such that $\frac{d^{+}l(t)}{dt}\leq
C\cdot l(t)$ for all $t\in(t_{1},t_{2})$. Once this is proved, then
$l(t)\leq0$ for all $t\in\lbrack t_{1},t_{2}]$ by Lemma
\ref{lem on dist is lower cont in Eucl} and Lemma \ref{lem on forward deriv}.
Hence $\tau(t)\in\mathcal{J}(t)\ for$ $t\in\lbrack t_{1},t_{2}]$, which is the
required contradiction.

Now our proof of the maximum principle for the time-dependent case diverges
from Hamilton's proof of the maximum principle for the time-independent case.
This is a necessity in our approach. The key difference is that we will not
use the general principle (Lemma \ref{lem taking deriv of sup funct}). We will
calculate $\frac{d^{+}l(t)}{dt}$ directly from the\textrm{ }definition. Also
our proof will not need the cutoff argument which appeared after Lemma
\ref{lem taking deriv of sup funct}. For any $t\in(t_{1},t_{2})$ there is a
sequence $s_{i}\rightarrow0^{+}$ such that%
\[
\frac{d^{+}l(t)}{dt}=\lim_{i\rightarrow\infty}\frac{l(t+s_{i})-l(t)}{s_{i}}.
\]

For\ any $v\in\partial\mathcal{J}(t)$, as in previous section we define
$S_{v}\subset\mathbb{R}^{n}$ to be the set of outward normal
directions\textrm{ }$n$\textrm{ }of the supporting hyperplanes of
$\mathcal{J}(t)$ at $v$;\textrm{ }we require that\textrm{ }$n$\textrm{ }be
unit with respect to the Euclidean metric. Define
\[
g(v,n,t)=n\cdot\lbrack\tau(t)-v].
\]
Since $\tau(t)\notin\mathcal{J}(t)$ for $t\in(t_{1},t_{2})$, we have
\[
l(t)=\sup_{v\in\partial\mathcal{J}(t)}\sup_{n\in S_{v}}g(v,n,t)
\]
and so we can find a sequence of points $v_{i}\in\partial\mathcal{J}(t+s_{i})$
and $n_{i}\in S_{v_{i}}$ such that $g(v_{i},n_{i},t+s_{i})=l(t+s_{i}%
)=|\tau(t+s_{i})-v_{i}|$. We can also find $v_{\infty}\in\partial
\mathcal{J}(t)$ and $n_{\infty}\in S_{v_{\infty}}$ such that $g(v_{\infty
},n_{\infty},t)=l(t)=|\tau(t)-v_{\infty}|$. It is not obvious that such
$v_{\infty}$ exists when $t=t_{1}$; this is one of the reason why we use Lemma
\ref{lem on forward deriv}. The proof below does not need a subsequence of
$v_{i}$ to converge to $v_{\infty}$ or a subsequence of $n_{i}$ to converge to
$n_{\infty}$.%
\begin{align*}
\frac{d^{+}l(t)}{dt}  &  =\lim_{i\rightarrow\infty}\frac{g(v_{i},n_{i}%
,t+s_{i})-g(v_{\infty},n_{\infty},t)}{s_{i}}\\
&  =\lim_{i\rightarrow\infty}\frac{n_{i}\cdot\lbrack\tau(t+s_{i}%
)-v_{i}]-n_{\infty}\cdot\lbrack\tau(t)-v_{\infty}]}{s_{i}}\\
&  =\lim_{i\rightarrow\infty}\frac{n_{i}\cdot\lbrack\tau(t+s_{i}%
)-\tau(t)]+n_{i}\cdot\tau(t)-n_{i}\cdot v_{i}-n_{\infty}\cdot\lbrack
\tau(t)-v_{\infty}]}{s_{i}}.
\end{align*}
Since $(F(v_{\infty},t),1)\in C_{(v_{\infty},t)}\mathcal{L}$, we can find a
subsequence $s_{i}$ and vectors $F_{i}\rightarrow F(v_{\infty},t)$ as
$i\rightarrow\infty$ such that $(v_{\infty}+s_{i}F_{i})\in\mathcal{J}%
(t+s_{i})$. Note that $v_{i}\in\partial\mathcal{J}(t+s_{i})$ and $n_{i}$ is
the outward normal direction of the supporting hyperplane at $v_{i}$. We have%
\[
n_{i}\cdot\lbrack v_{\infty}+s_{i}F_{i}-v_{i}]\leq0.
\]
Hence%
\begin{align*}
\frac{d^{+}l(t)}{dt}  &  =\lim_{i\rightarrow\infty}\,\,\{n_{i}\cdot
\lbrack\frac{\tau(t+s_{i})-\tau(t)}{s_{i}}-F_{i}]+\frac{n_{i}\cdot\lbrack
v_{\infty}+s_{i}F_{i}-v_{i}]}{s_{i}}\\
&  +\frac{(n_{i}-n_{\infty})\cdot\lbrack\tau(t)-v_{\infty}]}{s_{i}}\}\\
&  \leq\lim_{i\rightarrow\infty}\{n_{i}\cdot\lbrack\frac{\tau(t+s_{i}%
)-\tau(t)}{s_{i}}-F_{i}]+\frac{(n_{i}-n_{\infty})\cdot\lbrack\tau
(t)-v_{\infty}]}{s_{i}}\}\newline\\
&  \leq\lim_{i\rightarrow\infty}n_{i}\cdot\lbrack\frac{\tau(t+s_{i})-\tau
(t)}{s_{i}}-F_{i}]\newline\leq\lim_{i\rightarrow\infty}|\frac{\tau
(t+s_{i})-\tau(t)}{s_{i}}-F_{i}|\newline\\
&  =|F(\tau(t),t)-F(v_{\infty},t)|\newline\leq C\cdot|\tau(t)-v_{\infty
}|\newline=C\cdot l(t).
\end{align*}
We have used $(n_{i}-n_{\infty})\cdot\lbrack\tau(t)-v_{\infty}]\leq0$ to get
the second inequality above. This is because $n_{i}\cdot\lbrack\tau
(t)-v_{\infty}]\leq|\tau(t)-v_{\infty}|$ and $|\tau(t)-v_{\infty}|=n_{\infty
}\cdot\lbrack\tau(t)-v_{\infty}].$ We have used $|n_{i}|=1$ to get the third
inequality above. Proposition \ref{prop on ode presev time dep convex set} is
now proved.

The rest of this section is devoted to the proof of Theorem
\ref{Hamilton time-dep weak max thm}. We will prove it by contradiction.
Suppose we have a solution $\sigma(x,t)$ of the PDE (\ref{general heat eq 1.3}%
) on $[t_{0},T]$ which starts with $\sigma(x,t_{0})\in\mathcal{K}_{x}(t_{0})$
for all $x\in M$\textrm{ }and which goes out of space-time track $\mathcal{T}$
at some time $t_{2}$. Since $\mathcal{T}$ is closed, there is a time
$t_{1}\geq t_{0}$ such that $\sigma(x,t_{1})\in\mathcal{K}_{x}(t_{1})$ for all
$x\in M$ and for any $t_{1}<t<t_{2}$ there is $x$ such that $\sigma
(x,t)\notin\mathcal{K}_{x}(t)$. Below we will focus on the time interval
$[t_{1},t_{2}]$.

Define the function%
\begin{equation}
f(t)=\sup_{x\in M}\,d(\sigma(x,t),\mathcal{K}_{x}(t))\text{ for }t\in\lbrack
t_{1},t_{2}] \label{eq 2.6 for the dist in bundle}%
\end{equation}
where $d$ is distance on $V_{x}$ defined by the metric $h$. It is clear from
our choice that $f(t_{1})=0$, $f(t)>0$ for $t>t_{1}$. Note that $f(t)$ is not
necessarily continuous.

Next we prove a lemma which will enable us later to apply Lemma
\ref{lem on forward deriv} to $f(t)$ defined by
(\ref{eq 2.6 for the dist in bundle}). Let $\hat{\sigma}(x,t)$ be any
continuous section of bundle $V$\textrm{ }which satisfies that\textrm{ }%
$\hat{\sigma}(x,t_{1})\in\mathcal{K}_{x}(t_{1})$ for all\textrm{ }$x\in
M$\textrm{ }and where for each $t\in(t_{1},t_{2}]$ there is $x$ such that
$\hat{\sigma}(x,t)$ is not in $\mathcal{K}_{x}(t)$. We define the function
$\widehat{g}:M\times\lbrack t_{1},t_{2}]\rightarrow\mathbb{R}$ by%

\[
\widehat{g}(x,t)=d(\hat{\sigma}(x,t),\mathcal{K}_{x}(t)),
\]
and define the function%
\[
\hat{f}(t)=\sup_{x\in M}\widehat{g}(x,t)\text{ for }t\in\lbrack t_{1},t_{2}].
\]
By assumption $\hat{f}(t_{1})=0,$ and $\hat{f}(t)>0$ for $t\in(t_{1},t_{2}]$.
For any $t\in\lbrack t_{1},t_{2})$ and any sequence $s_{i}\rightarrow0^{+}$,
there is a subsequence $s_{i}$ and a sequence $x_{i}\in M$ such that
$\widehat{g}(x_{i},t+s_{i})=\sup_{x\in M}\,\widehat{g}(x,t+s_{i})$ and
$x_{i}\rightarrow x_{\infty}$.

\begin{lemma}
\label{lem on dist is lower cont in bundle} For the space-time track
$\mathcal{T}$ satisfying the assumption of Theorem
\ref{Hamilton time-dep weak max thm}, $\hat{f}(t)$ is left lower
semi-continuous and is right-continuous on $[t_{1},t_{2}]$, and for
$t\in\lbrack t_{1},t_{2})$\ the above chosen $x_{\infty}$ satisfies\textrm{ }
\[
\widehat{g}(x_{\infty},t)=\hat{f}(t).
\]

\end{lemma}

\begin{proof}
[Proof of the lemma.]First we show that $\widehat{f}(t)$ is lower
semi-continuous. $\widehat{f}(t)$ is obviously lower semi-continuous at
$t=t_{1}$. At any $t=t_{a}\in(t_{1},t_{2}]$, we have $\widehat{f}(t_{a})>0$.
We fix $x_{a}$ such that $\widehat{f}(t_{a})=\widehat{g}(x_{a},t_{a})$. Then
since $\mathcal{T}$ is closed, there is an $\varepsilon>0$ such that
$\hat{\sigma}(x_{a},t)\notin\mathcal{K}_{x_{a}}(t)$ for $t\in(t_{a}%
-\varepsilon,t_{a}+\varepsilon)$. We can apply Lemma
\ref{lem on dist is lower cont in Eucl} to $\widehat{g}(x_{a},t)$ in the fiber
$V_{x_{a}}$ to conclude that $\widehat{g}(x_{a},\cdot)$ is lower
semi-continuous at $t=t_{a}$. Hence for any $s_{i}\rightarrow0$%
\[
\lim\inf_{i\rightarrow+\infty}\widehat{f}(t_{a}+s_{i})\geq\lim\inf
_{i\rightarrow+\infty}\widehat{g}(x_{a},t_{a}+s_{i})\newline\geq\widehat
{g}(x_{a},t_{a})\newline=\widehat{f}(t_{a}).
\]
Hence $\widehat{f}(t)$ is lower semi-continuous at time $t=t_{a}$, and hence
on [$t_{1},t_{2}].$

To prove the right-continuity of $\widehat{f}(t)$, it suffices to prove the
upper right-continuity. For any $t_{a}\in\lbrack t_{1},t_{2})$ and any
sequence $s_{i}\rightarrow0^{+}$ we will show that there is a subsequence
$s_{i}$ such that $\lim_{i\rightarrow\infty}\widehat{f}(t_{a}+s_{i}%
)\leq\widehat{f}(t_{a})$. By passing to a subsequence if necessarily we may
assume that $\lim_{i\rightarrow\infty}\widehat{f}(t_{a}+s_{i})$ exists. Choose
$x_{i}\in M$ satisfying $\widehat{f}(t_{a}+s_{i})=\widehat{g}(x_{i}%
,t_{a}+s_{i})$; without loss of generality, we may assume that $x_{i}%
\rightarrow x_{\infty}$ by taking a subsequence if necessary. Let $v_{\infty
}\in\mathcal{K}_{x_{\infty}}(t_{a})$ such that $\widehat{g}(x_{\infty}%
,t_{a})=d(\widehat{\sigma}(x_{\infty},t_{a}),v_{\infty})$. If follows from
$\widehat{f}(t_{a}+s_{i})>0$, the invariance of $\mathcal{T}_{x}$ under
parallel translation and the closedness of $\mathcal{T}$, that $(v_{\infty
},t_{a})\in\partial\mathcal{T}_{x_{\infty}}$. By the assumption of Theorem
\ref{Hamilton time-dep weak max thm} and Proposition
\ref{prop on ode presev time dep convex set}, $C_{(v_{\infty},t_{a}%
)}\mathcal{T}_{x_{\infty}}$ is nonempty. Then there is a subsequence
$(v_{\infty}+s_{i}W_{i})\in\mathcal{K}_{x_{\infty}}(t_{a}+s_{i})$ with
$W_{i}\rightarrow W$ for some $W\in V_{x_{\infty}}$. Hence%
\begin{equation}
d(\widehat{\sigma}(x_{\infty},t_{a}+s_{i}),v_{\infty}+s_{i}W_{i})\geq
d(\widehat{\sigma}(x_{\infty},t_{a}+s_{i}),\mathcal{K}_{x_{\infty}}%
(t_{a}+s_{i})). \label{eq 2.7 temp use}%
\end{equation}
Since $\widehat{\sigma}(x,t_{a})$ is continuous in $x$ and $\mathcal{K}%
_{x}(t_{a}+s_{i})$ is invariant under parallel translation $\nabla(t_{a}%
+s_{i})$\ for any $x\in M$, $d(\widehat{\sigma}(x_{\infty},t_{a}%
+s_{i}),\mathcal{K}_{x_{\infty}}(t_{a}+s_{i}))$ can be chosen arbitrarily
close to $d(\widehat{\sigma}(x_{i},t_{a}+s_{i}),\mathcal{K}_{x_{i}}%
(t_{a}+s_{i}))=\widehat{g}(x_{i},t_{a}+s_{i})$ when $i$ is large, so the right
side of (\ref{eq 2.7 temp use}) approaches $\lim_{i\rightarrow+\infty}%
\widehat{f}(t_{a}+s_{i})$. The left side of (\ref{eq 2.7 temp use}) approaches%
\[
d(\widehat{\sigma}(x_{\infty},t_{a}),v_{\infty})=d(\widehat{\sigma}(x_{\infty
},t_{a}),\mathcal{K}_{x_{\infty}}(t_{a}))\leq\widehat{f}(t_{a}).
\]
Now we have proved $\lim_{i\rightarrow+\infty}\widehat{f}(t_{a}+s_{i}%
)\leq\widehat{f}(t_{a})$ and hence the right continuity of $\widehat{f}(t)$.

By taking the limit of (\ref{eq 2.7 temp use}) we have
\[
\widehat{g}(x_{\infty},t_{a})=d(\widehat{\sigma}(x_{\infty},t_{a}%
),\mathcal{K}_{x_{\infty}}(t_{a}))\geq\lim_{i\rightarrow+\infty}\widehat
{f}(t_{a}+s_{i}).
\]
Since $\widehat{f}(t)$ is right-continuous and $\widehat{g}(x_{\infty}%
,t_{a})\leq\widehat{f}(t_{a})$, we conclude that $\widehat{g}(x_{\infty}%
,t_{a})=\widehat{f}(t_{a})$ for any $t_{a}\in\lbrack t_{1},t_{2})$. The lemma
is proved.
\end{proof}

Now we go back to the proof of Theorem \ref{Hamilton time-dep weak max thm}.
Let $f(t)$ be the function defined in (\ref{eq 2.6 for the dist in bundle}),
we will prove that there is a constant $C<+\infty$ such that $\frac{d^{+}%
f(t)}{dt}\leq C\cdot f(t)$ for $t\in(t_{1},t_{2})$. Once this is proved, from
Lemma \ref{lem on dist is lower cont in bundle} and Lemma
\ref{lem on forward deriv} we conclude that $f(t)=0$ for $t\in\lbrack
t_{1},t_{2}],$ and hence $\sigma(x,t)\in\mathcal{K}_{x}(t)$ for all $x\in M$
and $t\in\lbrack t_{1},t_{2}]$. We get the required contradiction.

For any $t_{a}\in(t_{1},t_{2})$ there exists a sequence $s_{i}\rightarrow
0^{+}$ such that%
\[
\frac{d^{+}f(t_{a})}{dt}=\lim_{i\rightarrow\infty}\frac{f(t_{a}+s_{i}%
)-f(t_{a})}{s_{i}}.
\]

We define the function%
\[
g(x,v,n,t)=n\cdot\lbrack\sigma(x,t)-v],\text{ for }x\in M,n\in V_{x},\,v\in
V_{x},\text{ and }t\in\lbrack t_{1},t_{2}].
\]
For any $v\in\partial\mathcal{K}_{x}(t)$, we define $S_{v}\subset V_{x}$ to be
the set of the outward unit normal directions $n$ of the supporting
hyperplanes of $\mathcal{K}_{x}(t)$ in $V_{x}$ at $v$. Then, for any $t>t_{1}$
since $\mathcal{K}_{x}(t)$ is not empty and $\sigma(x,t)$ is not in the
interior of $\mathcal{K}_{x}(t)$ for some $x\in M$, it is well-known that%
\[
f(t)=\sup_{x\in M}\,\sup_{v\in\partial\mathcal{K}_{x}(t)}\,\sup_{n\in S_{v}%
}g(x,v,n,t).
\]
Note that the set over which we take the supremum in the definition of $f(t)$
depends on time. This is why we compute $\frac{d^{+}f(t_{a})}{dt}$ directly
rather than using Lemma \ref{lem taking deriv of sup funct}.

We can find a sequence of points $x_{i}\in M$, $v_{i}\in\partial
\mathcal{K}_{x_{i}}(t_{a}+s_{i})$, and $n_{i}\in S_{v_{i}}$, such that
$g(x_{i},v_{i},n_{i},t_{a}+s_{i})=f(t_{a}+s_{i})$, by Lemma
\ref{lem on dist is lower cont in bundle} we may assume $x_{i}\rightarrow
x_{\infty}\in M$ and $f(t_{a})=d(\sigma(x_{\infty},t_{a}),\mathcal{K}%
_{x_{\infty}}(t_{a}))$. Since $t_{a}>t_{1}$ and $f(t_{a})>0$, we have
$f(t_{a})=\sup_{v\in\partial\mathcal{K}_{x_{\infty}}(t_{a})}\,\sup_{n\in
S_{v}}g(x_{\infty},v,n,t_{a})$. Let $v_{\infty}\in\partial\mathcal{K}%
_{x_{\infty}}(t_{a})$ and $n_{\infty}\in S_{v_{\infty}}$, such that
$g(x_{\infty},v_{\infty},n_{\infty},t_{a})=f(t_{a})$.

We claim that there is a subsequence $i$ such that $v_{i}\rightarrow
v_{\infty}$ and $n_{i}\rightarrow n_{\infty}$ in the bundle $V$. Since
$d(\sigma(x_{i},t_{a}+s_{i}),v_{i})=f(t_{a}+s_{i})$ and $f(t_{a}+s_{i})$ is
uniformly bounded from above by the right-continuity of $f(t)$, we can rule
out the divergence of $v_{i}$ to $\infty$. We\textrm{ }may assume that there
is a subsequence $i$ such that $v_{i}\rightarrow\widehat{v}_{\infty}$ and
$n_{i}\rightarrow\widehat{n}_{\infty}$. By the closedness of $\mathcal{T}$ we
have $\hat{v}_{\infty}\in\mathcal{K}_{x_{\infty}}(t_{a})$, also we have
$|\widehat{n}_{\infty}|=1$. By taking the limit of%
\[
n_{i}\cdot\lbrack\sigma(x_{i},t_{a}+s_{i})-v_{i}]=f(t_{a}+s_{i}),\hspace
{0.25in}|\sigma(x_{i},t_{a}+s_{i})-v_{i}|=f(t_{a}+s_{i}),
\]
we get%
\begin{align}
\widehat{n}_{\infty}\cdot\lbrack\sigma(x_{\infty},t_{a})-\widehat{v}_{\infty
}]  &  =f(t_{a}),\label{eq 2.8 for temp use}\\
|\sigma(x_{\infty},t_{a})-\widehat{v}_{\infty}|  &  =f(t_{a}).
\label{eq 2.9 for temp use}%
\end{align}
By the convexity of $\mathcal{K}_{x_{\infty}}(t_{a})$,
(\ref{eq 2.9 for temp use}) implies $\widehat{v}_{\infty}=v_{\infty}$ and
(\ref{eq 2.8 for temp use}) implies $\widehat{n}_{\infty}=n_{\infty}$. The
claim is proved.

Then%
\begin{align*}
\frac{d^{+}f(t_{a})}{dt}  &  =\lim_{i\rightarrow\infty}\frac{g(x_{i}%
,v_{i},n_{i},t_{a}+s_{i})-g(x_{\infty},v_{\infty},n_{\infty},t_{a})}{s_{i}}\\
&  =\lim_{i\rightarrow\infty}\frac{n_{i}\cdot\lbrack\sigma(x_{i},t_{a}%
+s_{i})-v_{i}]-n_{\infty}\cdot\lbrack\sigma(x_{\infty},t_{a})-v_{\infty}%
]}{s_{i}}\\
&  =\lim_{i\rightarrow\infty}\left\{  \frac{n_{i}\cdot\lbrack\sigma
(x_{i},t_{a}+s_{i})-\sigma(x_{i},t_{a})]+n_{i}\cdot\sigma(x_{i},t_{a})}{s_{i}%
}\right. \\
&  +\left.  \frac{-n_{i}\cdot v_{i}-n_{\infty}\cdot\lbrack\sigma(x_{\infty
},t_{a})-v_{\infty}]}{s_{i}}\right\}  .
\end{align*}

To estimate quantities at $(x_{i},t_{a}+s_{i})$ and at $(x_{\infty},t_{a})$ in
$\frac{d^{+}f(t_{a})}{dt},$\textrm{ }we interpose quantities at $(x_{i}%
,t_{a})$ (see (\ref{eq for interpol}) below). By Lemma
\ref{lem on dist is lower cont in bundle}, $d(\sigma(x_{\infty},t_{a}%
),\mathcal{K}_{x_{\infty}}(t_{a}))=f(t_{a})>0,$ and hence $\sigma(x_{\infty
},t_{a})\notin\mathcal{K}_{x_{\infty}}(t_{a})$. It follows from $\mathcal{K}%
_{x}(t_{a})$ being invariant under parallel translation that for large enough
$i,$ $\sigma(x_{i},t_{a})\notin\mathcal{K}_{x_{i}}(t_{a})$. We can choose
$v_{i}^{\ast}\in\partial\mathcal{K}_{x_{i}}(t_{a})$ and $n_{i}^{\ast}\in
S_{v_{i}^{\ast}}$ such that $d(\sigma(x_{i},t_{a}),\mathcal{K}_{x_{i}}%
(t_{a}))=n_{i}^{\ast}\cdot\lbrack\sigma(x_{i},t_{a})-v_{i}^{\ast}]$. Such
$v_{i}^{\ast}$ and $n_{i}^{\ast}$ may not exist at time $t_{a}=t_{1}$ since
$\sigma(x_{\infty},t_{1})\in\mathcal{K}_{x_{\infty}}(t_{1})$; this is another
reason why we need Lemma \ref{lem on forward deriv}.

We claim that there is a sequence of vectors $F_{i}\in V_{x_{i}}$ such that
for any $\varepsilon>0$ there is an $i_{0}$ such that for any $i\geq i_{0}$ we
have $v_{i}^{\ast}+s_{i}F_{i}\in\mathcal{K}_{x_{i}}(t_{a}+s_{i})$ and
$|F_{i}-F(x_{i},v_{i}^{\ast},t_{a})|\leq\varepsilon$. The claim can be proved
by studying a family indexed by $i$ of ODE (\ref{general ode 1.4}) in
$V_{x_{i}}$ with initial time $t_{a}$ and initial value $\sigma_{x_{i}}%
(t_{a})=v_{i}^{\ast}$. We write the solution $\sigma_{x_{i}}(t_{a}%
+s_{i})=v_{i}^{\ast}+s_{i}F_{i}$. It follows from the assumption of Theorem
\ref{Hamilton time-dep weak max thm} that $\sigma_{x_{i}}(t_{a}+s_{i}%
)\in\mathcal{K}_{x_{i}}(t_{a}+s_{i}).$ Since $F(x,\sigma,t)$ is Lipschitz in
$\sigma$, the inequality $|F_{i}-F(x_{i},v_{i}^{\ast},t_{a})|\leq\varepsilon$
follows from the fact that solutions of ordinary differential equations depend
continuously on their parameters, in this case the parameters are $x_{i}\in M$
and $v_{i}^{\ast}\in\partial\mathcal{K}_{x_{i}}(t_{a})$ varying in compact domain.

Since $d(\sigma(x_{i},t_{a}),\mathcal{K}_{x_{i}}(t_{a}))=d(\sigma(x_{i}%
,t_{a}),v_{i}^{\ast})\leq f(t_{a})<\infty$, we can rule out the divergence of
$v_{i}^{\ast}$ to $\infty$. We may assume that a subsequence $v_{i}^{\ast}$
converges to $v_{\infty}^{\ast}\in V_{x_{\infty}}$, and get $d(\sigma
(x_{\infty},t_{a}),\mathcal{K}_{x_{\infty}}(t_{a}))=d(\sigma(x_{\infty}%
,t_{a}),v_{\infty}^{\ast}).$ By the closedness of the space-time track
$\mathcal{T}$ we have $v_{\infty}^{\ast}\in\mathcal{K}_{x_{\infty}}(t_{a})$.
Since
\[
d(\sigma(x_{\infty},t_{a}),v_{\infty}^{\ast})=d(\sigma(x_{\infty}%
,t_{a}),\mathcal{K}_{x_{\infty}}(t_{a}))=d(\sigma(x_{\infty},t_{a}),v_{\infty
})
\]
and $\mathcal{K}_{x_{\infty}}(t_{a})$ is convex, we conclude that $v_{\infty
}^{\ast}=v_{\infty}$. Our choice of $F_{i}$ ensures that $\lim_{i\rightarrow
\infty}F_{i}=F(x_{\infty},v_{\infty}^{\ast},t_{a})=F(x_{\infty},v_{\infty
},t_{a})$. Recall that $v_{i}\in\partial\mathcal{K}_{x_{i}}(t_{a}+s_{i})$ and
$n_{i}$ is the outward normal direction of the supporting hyperplane at
$v_{i}$. We have in each fiber $V_{x_{i}}$ and at time $t_{a}+s_{i}$%

\begin{equation}
n_{i}\cdot\lbrack v_{i}^{\ast}+s_{i}F_{i}-v_{i}]\leq0. \label{eq for interpol}%
\end{equation}

Hence%
\begin{align*}
\frac{d^{+}f(t_{a})}{dt}  &  =\lim_{i\rightarrow\infty}\,\{n_{i}\cdot
\lbrack\frac{\sigma(x_{i},t_{a}+s_{i})-\sigma(x_{i},t_{a})}{s_{i}}%
-F_{i}]+\frac{n_{i}\cdot\lbrack v_{i}^{\ast}+s_{i}F_{i}-v_{i}]}{s_{i}}\\
&  +\frac{n_{i}\cdot\lbrack\sigma(x_{i},t_{a})-v_{i}^{\ast}]-n_{\infty}%
\cdot\lbrack\sigma(x_{\infty},t_{a})-v_{\infty}]}{s_{i}}\}\\
&  \leq\lim_{i\rightarrow\infty}\{n_{i}\cdot\lbrack\frac{\sigma(x_{i}%
,t_{a}+s_{i})-\sigma(x_{i},t_{a})}{s_{i}}-F_{i}]\newline\\
&  +\frac{n_{i}\cdot\lbrack\sigma(x_{i},t_{a})-v_{i}^{\ast}]-n_{\infty}%
\cdot\lbrack\sigma(x_{\infty},t_{a})-v_{\infty}]}{s_{i}}\}\\
&  \leq\lim_{i\rightarrow\infty}\,\{n_{i}\cdot\lbrack\frac{\sigma(x_{i}%
,t_{a}+s_{i})-\sigma(x_{i},t_{a})}{s_{i}}-F_{i}]\},
\end{align*}
where to get the last inequality above we have used
\[
n_{i}\cdot\lbrack\sigma(x_{i},t_{a})-v_{i}^{\ast}]\leq n_{\infty}\cdot
\lbrack\sigma(x_{\infty},t_{a})-v_{\infty}].
\]
This is because
\[
n_{i}\cdot\lbrack\sigma(x_{i},t_{a})-v_{i}^{\ast}]\leq|\sigma(x_{i}%
,t_{a})-v_{i}^{\ast}|=d(\sigma(x_{i},t_{a}),\mathcal{K}_{x_{i}}(t_{a})),
\]
and at time $t_{a}$
\begin{align*}
d(\sigma(x_{i},t_{a}),\mathcal{K}_{x_{i}}(t_{a}))  &  \leq f(t_{a}%
)=d(\sigma(x_{\infty},t_{a}),\mathcal{K}_{x_{\infty}}(t_{a}))\\
&  =n_{\infty}\cdot\lbrack\sigma(x_{\infty},t_{a})-v_{\infty}]
\end{align*}
by our choice of $x_{\infty},\,v_{\infty}$, and $n_{\infty}$.%
\begin{align*}
\frac{d^{+}f(t_{a})}{dt}  &  \leq\lbrack\lim_{i\rightarrow\infty}n_{i}%
]\cdot\lbrack\lim_{i\rightarrow\infty}\frac{(\sigma(x_{i},t_{a}+s_{i}%
)-\sigma(x_{i},t_{a})}{s_{i}}-\lim_{i\rightarrow\infty}F_{i}]\\
&  =n_{\infty}\cdot\lbrack\frac{\partial}{\partial t}\sigma(x_{\infty}%
,t_{a})-F(x_{\infty},v_{\infty},t_{a})]\\
&  =n_{\infty}\cdot\lbrack\Delta(t_{a})\sigma(x_{\infty},t_{a})+u(x_{\infty
},t_{a})(\nabla(t_{a})\sigma(x_{\infty},t_{a}))\,\,\,\\
&  +F(x_{\infty},\sigma(x_{\infty},t_{a}),t_{a})-F(x_{\infty},v_{\infty}%
,t_{a})]\\
&  =n_{\infty}\cdot\lbrack\Delta(t_{a})\sigma(x_{\infty},t_{a})]+n_{\infty
}\cdot\lbrack u(x_{\infty},t_{a})(\nabla(t_{a})\sigma(x_{\infty}%
,t_{a}))]\newline\,\,\,\\
&  \,\,\,+n_{\infty}\cdot\lbrack F(x_{\infty},\sigma(x_{\infty},t_{a}%
),t_{a})-F(x_{\infty},v_{\infty},t_{a})].
\end{align*}

By the same argument as in section \ref{thm1,1 section} we conclude that
\begin{align*}
n_{\infty}\cdot\lbrack\Delta(t_{a})\sigma(x_{\infty},t_{a})]  &  \leq0,\\
n_{\infty}\cdot\lbrack u(x_{\infty},t_{a})(\nabla(t_{a})\sigma(x_{\infty
},t_{a}))]  &  =0.
\end{align*}
So%
\begin{align*}
\frac{d^{+}f(t_{a})}{dt}  &  \leq n_{\infty}\cdot\lbrack F(x_{\infty}%
,\sigma(x_{\infty},t_{a}),t_{a})-F(x_{\infty},v_{\infty},t_{a})]\\
&  \leq|F(x_{\infty},\sigma(x_{\infty},t_{a}),t_{a})-F(x_{\infty},v_{\infty
},t_{a})|\newline\\
&  \leq C\cdot|\sigma(x_{\infty},t_{a})-v_{\infty}|\newline=C\cdot f(t_{a}).
\end{align*}

Theorem \ref{Hamilton time-dep weak max thm} is proved.

\section{Proof of Theorem \ref{avoidance time-dep weak max thm}}

First we prove a version of Proposition
\ref{prop on ode presev time dep convex set} subject to an avoidance set.

\begin{proposition}
\label{prop on ode presev time dep convex avoid set} Let $U\subset
\mathbb{R}^{n}$\textrm{ }be an open subset, $\mathcal{J}(t)\subset
U,t\in\lbrack0,T]$\textrm{ }be a family of non-empty closed convex subsets and
$\mathcal{B}(t)\subset\mathcal{J}(t)$ be avoidance sets such that the
space-time track $\mathcal{L}$\ and the avoidance space-time track
$\mathcal{BL}=\{(v,t)\in\mathbb{R}^{n}\times\mathbb{R}:v\in\mathcal{B}%
(t),t\in\lbrack0,T]\}$\ are closed. Consider the ODE%
\begin{equation}
\frac{d\tau}{dt}=F(\tau,t), \label{eq 5.1 ode temp use}%
\end{equation}
where $F:U\times\lbrack0,T]\rightarrow\mathbb{R}^{n}$ is continuous in\textrm{
}$t$ and is Lipschitz in $\tau$. Then the following two statements are equivalent.

(i) For any $\,t_{0}\in\lbrack0,T)$ and any solution $\tau(t),t\in\lbrack
t_{0},T]$ of the ODE (\ref{eq 5.1 ode temp use}) with initial condition
$\tau(t_{0})\in\mathcal{J}(t_{0})\backslash\mathcal{B}(t_{0})$, either
$\tau(t)\in\mathcal{J}(t)$ for all $t\geq t_{0}$, or there is a time
$t_{1}>t_{0}$ such that $\tau(t)\in\mathcal{J}(t)\backslash\mathcal{B}(t)$ for
all $t\in\lbrack t_{0},t_{1})$ and $\tau(t_{1})\in\mathcal{B}(t_{1})$.

(ii) $(F(v,t),1)\in C_{(v,t)}\mathcal{L}$ for all $(v,t)\in(\partial
\mathcal{L})\backslash(\mathcal{BL})$.
\end{proposition}

\begin{proof}
This proposition can be proved as Proposition
\ref{prop on ode presev time dep convex set} except for the following issue
which arises in proving $(ii)\Longrightarrow(i)$. In the proof of Proposition
\ref{prop on ode presev time dep convex set} we have used the property
$(F(v,t),1)\in C_{(v,t)}\mathcal{L}$ for all $(v,t)\in\partial\mathcal{L}$,
however here this property holds only for $(v,t)\in(\partial\mathcal{L}%
)\backslash(\mathcal{BL})$. We need to ensure that $(v,t)$ can be chosen in
$(\partial\mathcal{L})\backslash(\mathcal{BL})$ when we use this property in
the proof of Proposition \ref{prop on ode presev time dep convex set}.

We adopt the notations used in the proof of Proposition
\ref{prop on ode presev time dep convex set} and resolve the issue. Since
$\mathcal{BL}$ is closed and the solution $\tau(t),t\in\lbrack t_{1},t_{2}]$
in the proof of Lemma \ref{lem on dist is lower cont in Eucl} does not enter
in $\mathcal{BL}$, there is a constant $\varepsilon>0$ such that
\[
\inf_{t\in\lbrack t_{1},t_{2}]}d(\tau(t),\mathcal{B}(t))\geq3\varepsilon.
\]
Since $(v_{t_{1}},t_{1})=(\tau(t_{1}),t_{1})\in(\partial\mathcal{L}%
)\backslash(\mathcal{BL}),l(t)$ is right-continuous at $t_{1}$ by the proof of
Lemma \ref{lem on dist is lower cont in Eucl}. Hence there is $t_{3}\in
(t_{1},t_{2})$ such that $f(t)\leq\varepsilon$ for all $t\in(t_{1},t_{3}).$
For any $t\in(t_{1},t_{3})$
\[
d(v_{t},\mathcal{B}(t))\geq d(\tau(t),\mathcal{B}(t))-d(v_{t},\tau
(t))\geq2\varepsilon,
\]
hence $(v_{t},t)\in(\partial\mathcal{L})\backslash(\mathcal{BL})$ for all
$t\in(t_{1},t_{3})$ and again $l(t)$ can be shown to be left lower
semi-continuous and right-continuous on $[t_{1},t_{3}]$.

For any $t\in(t_{1},t_{3}),$ choose the points $(v_{\infty},t)$ in
$\partial\mathcal{L}$ as in the proof of Proposition
\ref{prop on ode presev time dep convex set}. These points are at least
$2\varepsilon$ away from $\mathcal{BL}$,\textrm{ }so by statement (ii) we
still have the property $(F(v_{\infty},t),1)\in C_{(v_{\infty},t)}\mathcal{L}%
$, which was use in the proof of Proposition
\ref{prop on ode presev time dep convex set}. We may now repeat the rest of
the proof of Proposition \ref{prop on ode presev time dep convex set} to
conclude that there is a constant $C<+\infty$ such that $\frac{d^{+}f(t)}%
{dt}\leq C\cdot f(t)$ for all $t\in(t_{1},t_{3})$. By Lemma
\ref{lem on forward deriv} we get $l(t)=0$ on $[t_{1},t_{3}]$, which is the
required contradiction.
\end{proof}

The intuition behind the proof of Theorem
\ref{avoidance time-dep weak max thm} is as follows. Outside the avoidance set
(where the solution is assumed not to enter) the reaction term of the PDE
(i.e., corresponding to the associated ODE) wants to push the solution back
into the convex set. The diffusion part wants to keep the solution in the
convex set, possibly trying (but not succeeding) to push it into the avoidance part.

\begin{proof}
[Proof of Theorem \ref{avoidance time-dep weak max thm}]We will prove it by
contradiction. As in the proof of Theorem \ref{Hamilton time-dep weak max thm}%
, suppose we have a solution $\sigma(x,t)$ of PDE (\ref{general heat eq 1.3})
on $[t_{0},T]$ which starts with $\sigma(x,t_{0})\in\mathcal{K}_{x}%
(t_{0})\backslash\mathcal{A}_{x}(t_{0})$ for all\textrm{ }$x\in M$ and which
goes out of the space-time track\textrm{ }$\mathcal{T}$ at some time $t_{2}$.
Since $\mathcal{T}$ is closed, there is a time $t_{1}\geq t_{0}$ such that
$\sigma(x,t_{1})\in\mathcal{K}_{x}(t_{1})$ for all $x$ and for any $t\in
(t_{1},t_{2}]$ there is $x$ such that $\sigma(x,t)\notin\mathcal{K}_{x}(t)$.
Below we will focus on the time interval $[t_{1},t_{2}]$.

We define function%
\[
f(t)=\sup_{x\in M}d(\sigma(x,t),\mathcal{K}_{x}(t))\text{ for }t\in\lbrack
t_{1},t_{2}]
\]
where $d$ is the distance on $V_{x}$ defined by the metric $h$. It is clear
that $f(t_{1})=0$ and $f(t)>0$ for $t>t_{1}$.

Since the avoidance space-time track $\mathcal{AT}$\textrm{ }is closed
and$\ \sigma(x,t)\notin\mathcal{AT}$ for\textrm{ }all $x\in M$ and
$t\in\lbrack t_{1},t_{2}]$, there is an $\varepsilon>0$ such that%
\[
\inf_{x\in M,t\in\lbrack t_{1},t_{2}]}d(\sigma(x,t),\mathcal{A}_{x}%
(t))\geq3\varepsilon.
\]

By Proposition \ref{prop on ode presev time dep convex avoid set} we have
$(F(x,v,t),1)\in C_{(v,t)}(\mathcal{T}_{x})$ for all $(v,t)\in(\partial
\mathcal{T}_{x})\backslash(\mathcal{AT}_{x}),$ however we have used the
property $(F(x,v,t),1)\in C_{(v,t)}\mathcal{T}_{x}$ for all $(v,t)\in
\partial\mathcal{T}_{x}$ in the proof of Lemma
\ref{lem on dist is lower cont in bundle}, we need to modify the proof of
Lemma \ref{lem on dist is lower cont in bundle} to show that $f(t)$ is left
lower semi-continuous and right-continuous. We adopt the notations used in the
proof of Lemma \ref{lem on dist is lower cont in bundle} and replace
$\widehat{\sigma}(x,t)$ by $\sigma(x,t).$ When $t_{a}=t_{1},$ $(v_{\infty
},t_{1})=(\sigma(x,t_{1}),t_{1})\in(\partial\mathcal{T}_{x_{\infty}%
})\backslash(\mathcal{AT}_{x_{\infty}}),$ $f(t)$ is right-continuous at
$t_{1}$ by the same proof. Hence there is $t_{3}\in(t_{1},t_{2})$ such that
$f(t)\leq\varepsilon$ for all $t\in(t_{1},t_{3}).$ For any $t_{a}\in
(t_{1},t_{3})$
\[
d(v_{\infty},\mathcal{A}(t_{a}))\geq d(\sigma(x_{\infty},t_{a}),\mathcal{A}%
(t_{a}))-d(v_{\infty},\sigma(x_{\infty},t_{a}))\geq2\varepsilon,
\]
so $(v_{\infty},t_{a})\in(\partial\mathcal{T}_{x_{\infty}})\backslash
(\mathcal{AT}_{x_{\infty}})$ for all $t_{a}\in(t_{1},t_{3})$ and $f(t)$ is
left lower semi-continuous and right-continuous on $[t_{1},t_{3}]$.

We will prove that there is a constant $C<+\infty$\textrm{ }such that
$\frac{d^{+}f(t)}{dt}\leq C\cdot f(t)$ for all $t\in(t_{1},t_{3})$, then by
Lemma \ref{lem on forward deriv} we get $f(t)=0$ for all $t\in\lbrack
t_{1},t_{3}]$, which is the required contradiction.

For any $t\in(t_{1},t_{3})$ since $f(t)=\sup_{x\in M}d(\sigma(x,t),\mathcal{K}%
_{x}(t))<\varepsilon$, all the points in $\mathcal{T}$ we choose in the proof
of Theorem \ref{Hamilton time-dep weak max thm} are at least $2\varepsilon$
away from $\mathcal{AT}$,\textrm{ }so we can repeat the proof of Theorem
\ref{Hamilton time-dep weak max thm} to conclude that $\frac{d^{+}f(t)}%
{dt}\leq C\cdot f(t)$ for all $t\in(t_{1},t_{3})$. Hence Theorem
\ref{avoidance time-dep weak max thm} is proved.
\end{proof}

\end{document}